\documentclass{article}
\usepackage{amsthm}
\usepackage{amssymb}
\usepackage{amsxtra}
 
\newtheorem{thm}{Theorem}[section]
\newtheorem{lemma}{Lemma}[section]

\title{More precise pair correlation of zeros and primes in short intervals}
\author{Tsz Ho Chan}

\begin{document}
\maketitle
\begin{abstract}
Goldston and Montgomery [\ref{GM}] proved that the Strong Pair Correlation
Conjecture and two second moments of primes in short intervals are equivalent
to each other under Riemann Hypothesis. In this paper, we get the second main
terms for each of the above and show that they are almost equivalent to each
other.
\end{abstract}
\section{Introduction}
Riemann Hypothesis tells us that every non-trivial zero of the Riemann zeta
function is of the form $\rho = {1 \over 2} + i \gamma$. In the early 1970s, H.
Montgomery studied the distribution of the difference $\gamma - {\gamma}'$
between the zeros. Let
\begin{equation}
\label{basic}
F(X,T) = \mathop{\sum_{0 \leq \gamma \leq T}}_{0 \leq {\gamma}' \leq T}
         X^{i(\gamma - {\gamma}')} w(\gamma - {\gamma}')
\mbox{ and }
w(u) = {4 \over 4 + u^2}.
\end{equation}
Assuming Riemann Hypothesis, he proved in [\ref{M}] that, as $T \rightarrow
\infty$,
$$F(X,T) \sim {T \over 2 \pi} \log{X} + {T \over 2 \pi X^2} (\log{T})^2$$
for $1 \leq X \leq T$ (actually he only proved for $1 \leq X \leq o(T)$ and the
full range was done by Goldston [\ref{G}]). He conjectured that
$$F(X,T) \sim {T \over 2 \pi} \log{T}$$
for $T \leq X$ which is known as the Strong Pair Correlation Conjecture.

Later, Goldston and Montgomery [\ref{GM}] showed that, under Riemann Hypothesis,
the Strong Pair Correlation Conjecture and the following assertions are
equivalent to each other: As $X \rightarrow \infty$,
\begin{enumerate}
\item
$$\int_{1}^{X} \bigl(\psi(x+h)- \psi(x) - h \bigr)^2 dx \sim hX \log{{X \over
h}} \mbox{ for } X^{\epsilon} \leq h \leq X^{1 - \epsilon}.$$
\item
$$\int_{1}^{X} \bigl(\psi(x+ \delta x) - \psi(x) - \delta x \bigr)^2 dx \sim
{1 \over 2} \delta X^2 \log{1 \over \delta} \mbox{ for } X^{-1 + \epsilon} \leq
\delta \leq X^{-\epsilon}.$$
\end{enumerate}

Recently, Montgomery and Soundararajan [\ref{MS}] considered possible
cancellations among the error terms in Twin Prime Conjecture. They propose a
more precise asymptotic formula: For any $\epsilon > 0$,
$$\int_{1}^{X} \bigl(\psi(x+h)-\psi(x)-h \bigr)^2 dx = hX\log{X \over h} + BhX
+ o(hX)$$
for $X^{\epsilon} \leq h \leq X^{1-\epsilon}$, where $B = -C_{0}-\log{2\pi}$
and $C_{0}$ is Euler's constant. This more precise form suggests that one may
get a more precise form of the Strong Pair Correlation Conjecture by
employing the method of Goldston and Montgomery [\ref{GM}]. Loosely speaking,
our main result is that the following are ``equivalent'' under Riemann
Hypothesis:

(a) For every fixed $\epsilon > 0$,
$$\int_{1}^{X} \bigl(\psi(x+h)-\psi(x)-h\bigr)^2 dx = hX\log{X \over h} + BhX +
o(hX)$$
holds uniformly for $X^{\epsilon} \leq h \leq X^{1-\epsilon}$.

(b) For every fixed $\epsilon > 0$,
$$\int_{1}^{X} \bigl(\psi((1+\delta)x)-\psi(x)-\delta x\bigr)^2 dx =
   {1 \over 2}\delta X^2 \log{1 \over \delta} + C\delta X^2 + o(\delta X^2)$$
holds uniformly for $X^{-1+\epsilon} \leq \delta \leq X^{-\epsilon}$.

(c) For every fixed $\epsilon > 0$ and $A \geq 1+\epsilon$,
$$F(X,T) = {T \over 2 \pi}\log{T} + {D \over 2\pi}T +o(T)$$
holds uniformly for $T^{1+\epsilon} \leq X \leq T^A$ where $F$ is defined as 
(\ref{basic}).

Here $B=-C_0 - \log{2\pi}$, $C={1 \over 2}(B+1)$ and $D = -\log{2\pi} - 1$.
Hence, we deduce a more precise Strong Pair Correlation Conjecture
$$F(X,T) = {T \over 2\pi}\log{T \over 2\pi} - {T \over 2\pi} + o(T)$$
which has the same main terms as the number of zeros of the Riemann zeta
function up to height $T$ (which matches with the diagonal terms of $F(X,T)$).

As in [\ref{GM}], we shall use the following upper bounds by Saffari and Vaughan
[\ref{SV}]. Assume Riemann Hypothesis.
\begin{equation}
\label{ineq1}
\int_{1}^{X} \bigl(\psi(x + \delta x) - \psi(x) - \delta x \bigr)^2 dx \ll
\delta X^2 \Bigl(\log{2 \over \delta} \Bigr)^2
\end{equation}
for $0 < \delta \leq 1$, and
\begin{equation}
\label{ineq2}
\int_{1}^{X} \bigl(\psi(x+h) - \psi(x) - h \bigr)^2 dx \ll hX \Bigl(\log{2X 
\over h} \Bigr)^2
\end{equation}
for $0 < h \leq X$.
\section{(a) ``$\Leftrightarrow$'' (b)}
\begin{thm}
\label{thm1.1}
Assume Riemann Hypothesis. For every $0 < \epsilon < 1/2$, if
\begin{equation}
\label{a}
\int_{1}^{X} \bigl(\psi(x+h)-\psi(x)-h \bigr)^2 dx = hX\log{X \over h} + BhX+ 
O \Bigl({hX \over (\log{X/h})^5} \Bigr)
\end{equation}
holds uniformly for $X^{\epsilon} \leq h \leq X^{1-\epsilon}$, then for every
$0 < \epsilon < 1/2$,
\begin{equation}
\label{b}
\int_{1}^{X} \bigl(\psi(x+\delta x)-\psi(x)-\delta x \bigr)^2 dx =
{1 \over 2}\delta X^2 \log{1 \over \delta} + C\delta X^2 + 
O \Bigl({\delta X^2 \over \log{1/\delta}} \Bigr)
\end{equation}
holds uniformly for $X^{-1+\epsilon} \leq \delta \leq X^{-\epsilon}$.
Here $C=(1+B)/2$.
\end{thm}
Proof: The method is that of Saffari and Vaughan [\ref{SV}] employed in
Goldston and Montgomery [\ref{GM}]. For $\epsilon_1 > 0$, we want to deduce
from (\ref{a}) that 
\begin{equation}
\label{1.0}
\begin{split}
& \int_{0}^{\Delta} \int_{1}^{X} \bigl(\psi(x+\delta x)-\psi(x)-\delta x
\bigr)^2 dx d\delta \\
=& {1 \over 4}{\Delta}^2 X^2 \log{1 \over \Delta}+ K{\Delta}^2 X^2 + 
O \Bigl({{\Delta}^2 X^2 \over (\log{1/\Delta})^5} \Bigr)
\end{split}
\end{equation}
for $X^{-1+ \epsilon +\epsilon_1} \leq \Delta \leq X^{-\epsilon - \epsilon_1}$.
Here $K = 3/8 + B/4$. One can accomplish this by showing
\begin{equation}
\label{1.1}
\begin{split}
& \int_{0}^{\Delta} \int_{X/2}^{X} \bigl(\psi(x+\delta x)-\psi(x)-\delta x
\bigr)^2 dx d\delta \\
=& {3 \over 4} \Bigl({1 \over 4}{\Delta}^2 X^2 \log{1 \over 
\Delta}+ K {\Delta}^2 X^2 \Bigr) + O \Bigl({{\Delta}^2 X^2 \over 
(\log{1/ \Delta})^5} \Bigr)
\end{split}
\end{equation}
for $X^{-1 + \epsilon} \leq \Delta \leq X^{-\epsilon}$.
We substitute $h=\delta x$. Then the left hand side of (\ref{1.1}) becomes
\begin{equation}
\label{1.2}
\int_{\Delta X/2}^{\Delta X} \int_{h/\Delta}^{X} {f(x,h)^2 \over x} dx dh + 
\int_{X^{\epsilon / 2}}^{\Delta X/2} \int_{X/2}^{X} {f(x,h)^2 \over x} dx dh +
\int_{0}^{X^{\epsilon / 2}} \int_{X/2}^{X} {f(x,h)^2 \over x} dx dh
\end{equation}
where $f(x,h) = \psi(x+h)-\psi(x)-h$. By integration by parts, we have from 
(\ref{a}) that
\begin{eqnarray*}
\int_{U}^{V} {f(x,h)^2 \over x} dx &=& {h \over 2} \bigl((\log{V})^2 - 
(\log{U})^2 \bigr) + h(1+B)(\log{V}-\log{U}) \\
& &- h\log{h} (\log{V}-\log{U}) + O \Bigl({h \over (\log{X/h})^5} \Bigr)
\end{eqnarray*}
provided that $V^{\epsilon/2} \leq h \leq U^{1- {\epsilon/2}}$.
Here we need (\ref{a}) to hold uniformly for $X^{\epsilon} \leq h \leq
X^{1 - \epsilon}$ for every $\epsilon$.
Using this, the first term of (\ref{1.2}) becomes
$$\Bigl({\log{2} \over 8}-{3 \over 16}\Bigr){\Delta}^2 X^2 \log{\Delta} - 
{(\log{2})^2 \over 16}{\Delta}^2 X^2 - {(3+2B)\log{2} \over 16} {\Delta}^2 X^2$$
$$ + \Bigl({9 \over 32}+ {3B \over 16}\Bigr) {\Delta}^2 X^2 + 
O \Bigl({{\Delta}^2 X^2 \over (\log{1/ \Delta})^5} \Bigr)$$
if $X^{\epsilon/2} \leq h \leq ({h \over \Delta})^{1 -{\epsilon/2}}$
when ${\Delta X \over 2} \leq h \leq \Delta X$. One can easily check that this
is okay when $X^{-1+\epsilon} \leq \Delta \leq X^{-\epsilon}$.
The second term of (\ref{1.2}) becomes
$$-{\log{2} \over 8} {\Delta}^2 X^2 \log{\Delta} + {(\log{2})^2 \over 16} 
{\Delta}^2 X^2 + {(3+2B)\log{2} \over 16} {\Delta}^2 X^2 +
O \Bigl({{\Delta}^2 X^2 \over (\log{1/\Delta})^5} \Bigr)$$
if $X^{\epsilon/2} \leq h \leq ({X \over 2})^{1 -{\epsilon/2}}$ when
$X^{\epsilon/2} \leq h \leq {\Delta X \over 2}$. Again, one can check that
this is okay when $X^{-1+\epsilon} \leq \Delta \leq X^{-\epsilon}$.
By (\ref{ineq2}), the third term of (\ref{1.2}) is
$$\ll X^{\epsilon/2} X^{\epsilon} (\log{X})^2 \ll {{\Delta}^2 X^2 \over
(\log{1/\Delta})^5}.$$
Combining these, we have (\ref{1.1}) and hence (\ref{1.0}) by replacing $X$ by 
$X 2^{-k}$ in (\ref{1.1}), summing over $0 \leq k \leq M=[{7\log{\log{X}} \over
\log{2}}]$ and appealing to (\ref{ineq1}) with $X$ replaced by $X 2^{-M-1}$.
Note that when $X^{-1+\epsilon+\epsilon_1} \leq \Delta \leq X^{-\epsilon
-\epsilon_1}$, $(X 2^{-k})^{-1+\epsilon} \leq \Delta \leq (X 2^{-k})^{-\epsilon}
$ for $0 \leq k \leq M$. That is why we can use (\ref{1.1}) repeatedly.

We now deduce (\ref{b}) from (\ref{1.0}). For $X^{-1+\epsilon +2\epsilon_1} \leq
\Delta \leq X^{-\epsilon-2\epsilon_1}$, we set $\eta = {1 \over (\log{1/\Delta})
^4}$. By (\ref{1.0}) and Taylor's expansion for $\log{(1 + \eta)}$, we have
\begin{equation}
\label{e1}
\begin{split}
\int_{\Delta}^{(1 + \eta) \Delta} \int_{1}^{X} f(x,\delta x)^2 dx d\delta
&= {1 \over 2} \eta {\Delta}^2 X^2 \log{1 \over \Delta} + \Bigl(2K\eta - 
{1 \over 4} \eta \Bigr){\Delta}^2 X^2 \\
&+ O \Bigl({{\Delta}^2 X^2 \over (\log{1/\Delta})^5} \Bigr) +
O \Bigl({\eta}^2 {\Delta}^2 X^2 \log{1 \over \Delta} \Bigr).
\end{split}
\end{equation}
Let $g(x,\delta x) = f(x,\Delta x)$. By the identity
$$ f^2 - g^2 = 2f(f-g) - (f-g)^2 $$
and Cauchy's inequality, we find that
\begin{equation}
\label{e2}
\begin{split}
& \int_{\Delta}^{(1 + \eta) \Delta} \int_{1}^{X} f(x,\delta x)^2 - g(x,\delta x)
^2 dx d\delta \\
\ll& \Bigl(\int \int f(x,\delta x)^2 dx d\delta \Bigr)^{1/2}
\Bigl(\int \int \bigl(f(x,\delta x)-g(x,\delta x) \bigr)^2 dx d\delta \Bigr)^
{1/2} \\
&+ \int \int \bigl(f(x,\delta x)-g(x,\delta x) \bigr)^2 dx d\delta.
\end{split}
\end{equation}
But $f(x,\delta x) - g(x,\delta x) = f((1+\Delta)x,(\delta-\Delta)x)$, so
\begin{equation}
\label{e3}
\begin{split}
& \int_{\Delta}^{(1 + \eta) \Delta} \int_{1}^{X}
\bigl(f(x,\delta x)-g(x,\delta x)\bigr)^2 dx d\delta \\
=& \int_{0}^{\eta \Delta \over 1+\Delta} \int_{1+\Delta}^{(1+\Delta)X}
f(x,\delta x)^2 dx d\delta \ll {\eta}^2 {\Delta}^2 X^2 \log{1 \over \Delta \eta}
\end{split}
\end{equation}
by (\ref{1.0}), our choice for $\eta$ and the range for $\Delta$. Hence by
(\ref{e1}), (\ref{e2}) and (\ref{e3}), we have
\begin{eqnarray*}
& &\eta \Delta \int_{1}^{X} \bigl(\psi(x+\Delta x) - \psi(x) - \Delta x \bigr)^
2 dx =
\int_{\Delta}^{(1 + \eta) \Delta} \int_{1}^{X} g(x,\delta x)^2 dx d\delta \\
&=& \int_{\Delta}^{(1 + \eta) \Delta} \int_{1}^{X} f(x,\delta x)^2 dx d\delta +
O \Bigl({\eta}^{3/2} {\Delta}^2 X^2 \log{1 \over \Delta}\Bigr) \\
&=& {1 \over 2} \eta {\Delta}^2 X^2 \log{1 \over \Delta}+
\Bigl(2K-{1 \over 4} \Bigr) \eta {\Delta}^2 X^2 + 
O \Bigl({\eta}^{3/2} {\Delta}^2 X^2 \log{1 \over \Delta} \Bigr) + 
O \Bigl({{\Delta}^2 X^2 \over (\log{1/\Delta})^5} \Bigr).
\end{eqnarray*}
We divide both sides by $\eta \Delta$, and obtain
$$\int_{1}^{X} \bigl(\psi(x+\Delta x) - \psi(x) - \Delta x \bigr)^2 dx =
{1 \over 2}\Delta X^2 \log{1 \over \Delta} + \Bigl(2K-{1 \over 4}\Bigr)\Delta X^2
$$
$$+ O \Bigl({\eta}^{1/2} \Delta X^2 \log{1 \over \Delta}\Bigr) +
O \Bigl({\Delta X^2 \over \eta (\log{1/\Delta})^5} \Bigr).$$
Recall $\eta = {1 \over \log{(1/\Delta)}^4}$, we have, for $X^{-1+\epsilon
+2\epsilon_1} \leq \Delta \leq X^{-\epsilon-2\epsilon_1}$,
$$\int_{1}^{X} \bigl(\psi((1+\Delta)x) - \psi(x) - \Delta x \bigr)^2 dx =
{1 \over 2}\Delta X^2 \log{1 \over \Delta} + C\Delta X^2 +
O \Bigl({\Delta X^2 \over \log{1/\Delta}} \Bigr)$$ 
as $2K-{1 \over 4} = {1+B \over 2} = C$. Since, $\epsilon_1$ is arbitrary, we
have the theorem.

\smallskip

One can carry out a similar calculation from (\ref{b}) to (\ref{a}) (see
[\ref{GM}] for outline). Then one obtains the following converse:
\begin{thm}
\label{thm1.2}
Assume Riemann Hypothesis. For every $0 < \epsilon < 1/2$, if
\begin{equation}
\label{c}
\int_{1}^{X} \bigl(\psi(x+\delta x)-\psi(x)-\delta x \bigr)^2 dx =
   {1 \over 2}\delta X^2 \log{1 \over \delta} + C\delta X^2 + 
   O \Bigl({\delta X^2 \over (\log{1/\delta})^5} \Bigr)
\end{equation}
holds uniformly for $X^{-1+\epsilon} \leq \delta \leq X^{-\epsilon}$, then for
every $0< \epsilon < 1/2$,
\begin{equation}
\label{d}
\int_{1}^{X} \bigl(\psi(x+h)-\psi(x)-h \bigr)^2 dx = hX\log{X \over h} + BhX + 
O \Bigl({hX \over \log{X/h}} \Bigr)
\end{equation}
holds uniformly for $X^{\epsilon} \leq h \leq X^{1-\epsilon}$. Again 
$C=(1+B)/2$.
\end{thm}

The error terms for (\ref{a}), (\ref{b}), (\ref{c}) and (\ref{d}) are
reasonable because, for instance, from Montgomery and Soundararajan [\ref{MS}],
one expects that the error term for (\ref{a}) and (\ref{d}) to be $O(h
X^{1-{\epsilon}'})$. Similarly, one would expect the error term for (\ref{b})
and (\ref{c}) to be $O(\delta X^{2-{\epsilon}'})$. Thus, it would be nice if we
can have variations of Theorem \ref{thm1.1} and Theorem \ref{thm1.2} with these
(smaller) error terms. But, in order to do that, one has to be more careful
about the argument and one has to shorten the range for $\delta$ and $h$ a 
little bit, respectively.

\begin{thm}
Assume Riemann Hypothesis. For every $0 < \epsilon < 1/2$, if, for some
${\epsilon}_1 > 0$ ($\epsilon_1$ may depend on $\epsilon$),
\begin{equation}
\label{e}
\int_{1}^{X} \bigl(\psi(x+h)-\psi(x)-h \bigr)^2 dx = hX\log{X \over h} + BhX + 
O(hX^{1-{\epsilon}_1})
\end{equation}
holds uniformly for $X^{\epsilon} \leq h \leq X^{1-\epsilon}$, then
\begin{equation}
\label{f}
\int_{1}^{X} \bigl(\psi(x+\delta x)-\psi(x)-\delta x \bigr)^2 dx =
   {1 \over 2}\delta X^2 \log{1 \over \delta} + C\delta X^2 + 
   O(\delta X^{2-{\epsilon}_2})
\end{equation}
holds uniformly for $X^{-1+\epsilon+2\epsilon_1} \leq \delta \leq 
X^{-\epsilon}/2$ for some ${\epsilon}_2 > 0$. Here $C=(1+B)/2$.
\end{thm}
Proof: Let $f(x,h) = \psi(x+h)-\psi(x)-h$. Let $X^{-1+\epsilon+{\epsilon}_1}
\leq \Delta \leq X^{-\epsilon}$. By integration by parts, we have from 
(\ref{e}) that $$\int_{U}^{V} {f(x,h)^2 \over x} dx = {h \over 2} \bigl(
(\log{V})^2  - (\log{U})^2 \bigr) + h(1+B)(\log{V}-\log{U})$$
$$- h\log{h} (\log{V}-\log{U}) + O(h U^{-{\epsilon}_1})$$
as long as $V^\epsilon \leq h \leq U^{1-\epsilon}$ when $U \leq V \leq 2U$.
Then, by substituting $h = \delta x$,
\begin{eqnarray*}
& &\int_{0}^{\Delta} \int_{V/2}^{V} f(x,\delta x)^2 dx d\delta \\
&=&\int_{\Delta V/2}^{\Delta V} \int_{h/\Delta}^{V} {f(x,h)^2 \over x} dx dh + 
\int_{V^\epsilon}^{\Delta V/2} \int_{V/2}^{V} {f(x,h)^2 \over x} dx dh +
\int_{0}^{V^\epsilon} \int_{V/2}^{V} {f(x,h)^2 \over x} dx dh \\
&=&{3 \over 4} \Bigl({1 \over 4}{\Delta}^2 V^2 \log{1 \over \Delta} +
\bigl({3 \over 8} + {B \over 4} \bigr) {\Delta}^2 V^2 \Bigr) + 
O({\Delta}^2 V^{2-{\epsilon}_1}) + O \bigl(V^{2\epsilon}(\log{2V})^2 \bigr)
\end{eqnarray*}
where the last error term comes from evaluating the second integral at
$V^{\epsilon}$ and estimating the third integral by (\ref{ineq2}). The above is
okay as long as $V^{\epsilon} \leq h \leq ({h \over \Delta})^{1-\epsilon}$ when
${\Delta V \over 2} \leq h \leq \Delta V$, and $V^\epsilon \leq h \leq ({V \over
 2})^{1-\epsilon}$ when $V^{\epsilon} \leq h \leq {\Delta V \over 2}$. 

Now, suppose $X^{1-.6{\epsilon}_1} \leq V \leq X$. One has
$ V^\epsilon \leq X^\epsilon \leq {X^{-1+\epsilon+{\epsilon}_1}
X^{1-.6{\epsilon}_1} \over 2} \leq {\Delta V \over 2} \leq {V^{-\epsilon} V
\over 2} \leq \Bigl({V \over 2} \Bigr)^{1-\epsilon}$.
One only needs to check $h \leq ({h \over \Delta})^{1-\epsilon}$ when ${\Delta 
V \over 2} \leq h \leq \Delta V$. But $\Delta^{1-\epsilon}
h^{\epsilon} \leq \Delta^{1-\epsilon} (\Delta V)^{\epsilon} \leq \Delta
X^{\epsilon} \leq 1$. Thus
$$\int_{V/2}^{V} \int_{0}^{\Delta} \bigl(\psi(x+\delta x)-\psi(x)-\delta x
\bigr)^2
d\delta dx = {3 \over 4} \Bigl({1 \over 4}{\Delta}^2 V^2 \log{1 \over \Delta} +
K {\Delta}^2 V^2 \Bigr) + O({\Delta}^2 V^{2-{\epsilon}_1})$$
where $K = {3 \over 8} + {B \over 4}$ as long as $X^{1-.6{\epsilon}_1} \leq
V \leq X$. Replacing $V$ by $X 2^{-k}$ in above, summing over $0 \leq k \leq M 
= [{.6{\epsilon}_1 \log{X} \over \log{2}}]$ and appealing to (\ref{ineq1}) with 
$V = X 2^{-M-1}$,
$$\int_{0}^{\Delta} \int_{1}^{X} \bigl(\psi(x+\delta x)-\psi(x)-\delta x \bigr)
^2 dx
d\delta = {1 \over 4}{\Delta}^2 X^2 \log{1 \over \Delta} + K{\Delta}^2 X^2 +
O({\Delta}^2 X^{2-{\epsilon}_1}).$$
Consider $X^{-1+\epsilon+ 2\epsilon_1} \leq \Delta \leq X^{-\epsilon}/2$ and
set $\eta = X^{-2\epsilon_1 /3}$.
One just mimics the reduction in Theorem \ref{thm1.1} and obtain
\begin{eqnarray*}
\int_{1}^{X} \bigl(\psi(x+\Delta x) - \psi(x) - \Delta x \bigr)^2 dx &=&
{1 \over 2}\Delta X^2 \log{1 \over \Delta} + \Bigl(2K-{1 \over 4}\Bigr)
\Delta X^2 \\
& &+ O \Bigl({\eta}^{1/2} \Delta X^2 \log{1 \over \Delta} \Bigr) +
O \Bigl({\Delta X^{2-{\epsilon}_1} \over \eta} \Bigr) \\
&=& {1 \over 2}\Delta X^2 \log{1 \over \Delta} + C \Delta X^2 +
O(\Delta X^{2-{\epsilon}_2})
\end{eqnarray*}
for $X^{-1+\epsilon+2\epsilon_1} \leq \Delta \leq X^{-\epsilon}/2$, whenever
$0 < {\epsilon}_2 < {\epsilon}_1 / 3$.

\bigskip
Conversely, one has
\begin{thm}
Assume Riemann Hypothesis. For every $0 < \epsilon < 1/2$, if, for some
${\epsilon}_1 > 0$,
\begin{equation}
\int_{1}^{X} \bigl(\psi(x+\delta x)-\psi(x)-\delta x \bigr)^2 dx =
   {1 \over 2}\delta X^2 \log{1 \over \delta} + C\delta X^2 + 
   O(\delta X^{2-{\epsilon}_1})
\end{equation}
holds uniformly for $X^{-1+\epsilon} \leq \delta \leq X^{-\epsilon}$, then
\begin{equation}
\int_{1}^{X} \bigl(\psi(x+h)-\psi(x)-h \bigr)^2 dx = hX\log{X \over h} + BhX + 
O(hX^{1-{\epsilon}_2})
\end{equation}
holds uniformly for $X^{\epsilon+\epsilon_1} \leq h \leq X^{1-\epsilon-
\epsilon_1}$ for some ${\epsilon}_2 > 0$. Here $C=(1+B)/2$.
\end{thm}
\section{Preparation for (b) ``$\Leftrightarrow$'' (c)}
We need some Abelian-Tauberian results with two main terms and explicit error
terms. Basically, they are modified versions of the lemmas in Goldston and
Montgomery [\ref{GM}].
\begin{lemma}
\label{lemma2.1}
If $I(Y) = \int_{-\infty}^{\infty} e^{-2|y|} f(Y+y) dy = 1 + O({1 \over Y^8})$,
and if $f(y) \geq 0$ for all y as well as $\int_{0}^{1} f(Y+y) dy \ll 1$, then
\begin{equation}
\label{l1}
\int_{0}^{\log{2}} e^{2y}f(Y+y) dy = {3 \over 2} + O \Bigl({1 \over Y^2} \Bigr).
\end{equation}
\end{lemma}

Proof: Let $K_c(y) = max(0, c - |y|)$. We would choose $c \approx {1 \over
Y^2}$. From Lemma 1 of Goldston and Montgomery [\ref{GM}], one has
$$K_c(y) = {1 \over 2}e^{-2|y|} - {1 \over 4}e^{-2|y-c|} - {1 \over 4} 
e^{-2|y+c|} + \int_{-c}^{c} (c - |z|) e^{-2|y-z|} dz.$$
For $0 \leq a \leq 1$, let $L_{a,c}(y) = K_c(y-a)$, then one has
\begin{eqnarray*}
\int_{-\infty}^{\infty} L_{a,c}(y) f(Y+y) dy &=& {1 \over 2}I(Y+a) - {1 \over
4}I(Y+a+c) - {1 \over 4}I(Y+a-c) \\
& &+ \int_{-c}^{c} (c-|z|)I(Y+a+z) dz \\
&=& c^2 + O \Bigl({1 \over Y^8} \Bigr).
\end{eqnarray*}
Since,
$${1 \over \eta} \bigl(L_{a,c}(y) - L_{a,c - \eta}(y) \bigr) \leq 
\chi_{[a-c,a+c]}(y) \leq
{1 \over \eta} \bigl(L_{a,c+\eta}(y) - L_{a,c}(y) \bigr)$$
and $f \geq 0$,
$$\int_{-\infty}^{\infty} \chi_{[a-c,a+c]}(y) f(Y+y) dy = 2c + O(\eta) + O
\Bigl({1 \over \eta Y^8} \Bigr).$$
Taking $\eta = {1 \over Y^4}$, we have
$$\int_{-\infty}^{\infty} \chi_{[a-c,a+c]}(y) f(Y+y) dy =
\int_{-\infty}^{\infty} \chi_{[a-c,a+c]}(y) dy + O\Bigl({1 \over Y^4}\Bigr).$$
Now, we approximate $e^{2y}$ over $0 \leq y \leq \log{2}$ by a sum of $N=[Y^2]$
step functions, called the sum $S(y)$, each support on $[{n\log{2} \over N},
{(n+1)\log{2} \over N}]$ with value $e^{(2n\log{2}/N)}$ for $n = 0,1,...,N-1$.
Then, since $f \geq 0$ and $\int_{0}^{1} f(Y+y) dy \ll 1$,
$$\int_{0}^{\log{2}} e^{2y} f(Y+y) dy - \int_{0}^{\log{2}} S(y) f(Y+y) dy = 
O\Bigl({1 \over Y^2}\Bigr),$$
$$\int_{0}^{\log{2}} e^{2y} dy - \int_{0}^{\log{2}} S(y) dy = O\Bigl({1 \over
Y^2}\Bigr).$$
So,
\begin{eqnarray*}
\int_{0}^{\log{2}} e^{2y} f(Y+y) dy &=&\int_{-\infty}^{\infty} S(y) f(Y+y) dy +
O \Bigl({1 \over Y^2} \Bigr) \\
&=&\int_{-\infty}^{\infty} S(y) dy + O \Bigl({N \over Y^4} \Bigr) + O \Bigl({1 
\over Y^2} \Bigr) \\
&=&\int_{0}^{\log{2}} e^{2y} dy + O \Bigl({1 \over Y^2} \Bigr)
\end{eqnarray*}
which gives the lemma.
\begin{lemma}
\label{lemma2.2}
Let $f(t)$ be a continuous non-negative function defined for all $t\geq 0$ with
$f(t) \ll (\log{(t+2)})^2$. If $J(T) = \int_{0}^{T} f(t) dt =
T\log{T} + DT + \epsilon(T) T$ where $\epsilon(T) = O({1 \over (\log{T})
^\lambda})$ for ${\kappa}^{-1}(\log{1 \over \kappa})^{-(\lambda+2)} \leq T \leq
{\kappa}^{-1} (\log{1 \over \kappa})^{(\lambda+2)}$, with some $\lambda > 0$,
then
\begin{equation}
\label{l2}
\int_{0}^{\infty} \Bigl({\sin{\kappa u} \over u}\Bigr)^2 f(u) du = {\pi \over 
2} \kappa \log{1 \over \kappa} + \pi C' \kappa + O_{\lambda} \Bigl({\kappa 
\log{\log{1/\kappa}} \over (\log{1 / \kappa})^\lambda} \Bigr).
\end{equation}
Here $C'$ and $D$ are related by 
$C' = {1 \over 2}D - {1 \over 2}(C_0+\log{2})+1$.
\end{lemma}

Proof: Just like Lemma 2 of Goldston and Montgomery [\ref{GM}], we divide the
range of integration of (\ref{l2}) into four subintervals:
\begin{eqnarray*}
I_1 &=& \{ u : 0 \leq u \leq U_1={\kappa}^{-1} \Bigl(\log{1 \over
\kappa} \Bigr)^{-(\lambda+2)} \} \\
I_2 &=& \{ u : U_1 \leq u \leq U_2=c{\kappa}^{-1} \} \\
I_3 &=& \{ u : U_2 \leq u \leq U_3={\kappa}^{-1} \Bigl(\log{1 \over
\kappa} \Bigr)^{(\lambda+2)} \} \\
I_4 &=& \{ u : U_3 \leq u \leq \infty \}
\end{eqnarray*}
Here $c$, $(\log{1 \over \kappa})^{-(\lambda+2)} \leq c \leq (\log{1 \over 
\kappa})^{(\lambda+2)}$, is a parameter to be chosen later.

Since $f(t) \ll (\log{(t+2)})^2$, we see that
\begin{equation}
\label{2.2.1}
\int_{I_1} \ll \int_{0}^{U_1} {\kappa}^2 (\log{(u+2)})^2 du \ll {\kappa}^2 U_1
(\log{U_1})^2 \ll {\kappa \over (\log{1 \over \kappa})^{\lambda}},
\end{equation}
\begin{equation}
\label{2.2.2}
\int_{I_4} \ll \int_{U_3}^{\infty} u^{-2} (\log{u})^2 du \ll {U_3}^{-1}
(\log{U_3})^2 \ll {\kappa \over (\log{1 \over \kappa})^{\lambda}}.
\end{equation}
By writing $f(u) = \log{1 \over \kappa}+\log{\kappa u}+(D+1) + (f(u) -\log{u} -
(D+1))$, we can express the integral on $I_2$ as a sum of four integrals. Note
that
\begin{eqnarray*}
\int_{U_1}^{U_2} \Bigl({\sin{\kappa u} \over u} \Bigr)^2 du &=& \int_{0}^
{\infty} \Bigl({\sin{\kappa u} \over u}\Bigr)^2 du + O \Bigl({\kappa \over 
(\log{1 \over \kappa})^{(\lambda+2)}} \Bigr) + O \Bigl({\kappa \over c} \Bigr)
\\
&=& {\pi \over 2}\kappa + O \Bigl({\kappa \over (\log{1 \over \kappa})^
{(\lambda+2)}} \Bigr) +  O \Bigl({\kappa \over c}\Bigr),
\end{eqnarray*}
\begin{eqnarray*}
\int_{U_1}^{U_2} \Bigl({\sin{\kappa u} \over u}\Bigr)^2 \log{\kappa u} du &=&
\kappa \int_{(\log{1 / \kappa})^{-(\lambda+2)}}^{c} \Bigl({\sin{v} \over v}
\Bigr)^2 \log{v} dv \\ &=&
\kappa \int_{0}^{\infty} \Bigl({\sin{v} \over v}\Bigr)^2 \log{v} dv + O \Bigl(
{\kappa \log{\log{1 \over \kappa}} \over (\log{1 \over \kappa})^{(\lambda+2)}} 
\Bigr) + O\Bigl({\kappa \log{c} \over c}\Bigr).
\end{eqnarray*}
Put $r(u) = J(u) - u\log{u} - Du$, then $r'(u) = f(u) - \log{u} -(D+1)$. By
integrating by parts, we have
\begin{eqnarray*}
& &\int_{U_1}^{U_2} \Bigl({\sin{\kappa u} \over u}\Bigr)^2 \bigl(f(u)- \log{u}
-(D+1) \bigr) du = 
\int_{U_1}^{U_2} \Bigl({\sin{\kappa u} \over u}\Bigr)^2 dr(u) \\
&=& \Bigl[({\sin{\kappa u} \over u})^2 r(u)\Bigr]_{U_1}^{U_2} -\int_{U_1}^{U_2}
r(u) d\Bigl({\sin{\kappa u} \over u}\Bigr)^2 \\
&=& O\Bigl({\kappa \over (\log{1 \over \kappa})^{\lambda}}\Bigr) + 
O \Bigl(\mathop{max}_{U_1 \leq u \leq U_2}
|\epsilon(u)| \Bigl[\int_{U_1}^{U_2} |{\kappa \sin{\kappa u} \cos{\kappa u}
\over u}| du + \int_{U_1}^{U_2} \bigl({\sin{\kappa u} \over u}\bigr)^2 du 
\Bigr]\Bigr)
\end{eqnarray*}
with $c > 1$. Estimating the remaining integrals and using $\epsilon(u) 
\ll {1 \over (\log{u})^{\lambda}}$, we have
$$\int_{U_1}^{U_2} \Bigl({\sin{\kappa u} \over u}\Bigr)^2 \bigl(f(u)-\log{u}
-(D+1) \bigr) du \ll_{\lambda}
{\kappa \log{c} \over (\log{1 \over \kappa})^{\lambda}} + {\kappa \log{\log{1
\over \kappa}} \over (\log{1 \over \kappa})^{\lambda}}.$$
Hence,
\begin{equation}
\label{2.2.3}
\begin{split}
\int_{I_2} =& {\pi \over 2}\kappa \log{1 \over \kappa} + \kappa 
\int_{0}^{\infty} \Bigl({\sin{v} \over v}\Bigr)^2 \log{v} dv + {\pi \over 2}
(D+1) \kappa \\
&+ O_{\lambda} \Bigl({\kappa \log{1 \over \kappa} \over c} + {\kappa \log{c} 
\over c} + {\kappa \log{c} \over (\log{1 \over \kappa})^{\lambda}} + {\kappa 
\log{\log{1 \over \kappa}} \over (\log{1 \over \kappa})^{\lambda}}\Bigr).
\end{split}
\end{equation}
Over $I_3$, we write $f(u) = \log{1 \over \kappa} + \log{\kappa u} + (D+1)
 + (f(u) - \log{u} - (D+1))$ and break down the integral into four pieces
again. Note that
$$\int_{U_2}^{U_3} \Bigl({\sin{\kappa u} \over u}\Bigr)^2 du = \kappa \int_{c}^
{(\log{1 / \kappa})^{(\lambda+2)}} \Bigl({\sin{v} \over v}\Bigr)^2 dv \ll 
{\kappa \over c},$$
$$\int_{U_2}^{U_3} \Bigl({\sin{\kappa u} \over u}\Bigr)^2 \log{\kappa u} du = 
\kappa \int_{c}^{(\log{1 / \kappa})^{(\lambda+2)}} \Bigl({\sin{v} \over v}
\Bigr)^2 \log{v} dv \ll {\kappa \log{c} \over c},$$
and by a similar calculation as $I_2$,
$$\int_{U_2}^{U_3} \Bigl({\sin{\kappa u} \over u}\Bigr)^2
\bigl(f(u)-\log{u}-(D+1)\bigr) du \ll_{\lambda}
{\kappa \log{\log{1 \over \kappa}} \over (\log{1 \over \kappa})^{\lambda}}.$$
Thus,
\begin{equation}
\label{2.2.4}
\int_{I_3} \ll_{\lambda} {\kappa \log{1 \over \kappa} \over c} + {\kappa
\log{c} \over c} + {\kappa \log{\log{1 \over \kappa}} \over (\log{1 \over
\kappa})^{\lambda}}.
\end{equation}
Therefore, choosing $c \approx (log{1 \over \kappa})^{(\lambda+1)}$, we have 
from (\ref{2.2.1}), (\ref{2.2.2}), (\ref{2.2.3}) and (\ref{2.2.4}) that
$$\int_{0}^{\infty} \Bigl({\sin{\kappa u} \over u}\Bigr)^2 f(u) du = {\pi \over
2} \kappa \log{1 \over \kappa} + \Bigl[{\pi \over 2}(D+1) - {\pi \over 2}(C_0 +
\log{2}  - 1) \Bigr] \kappa + O_{\lambda} \Bigl({\kappa \log{\log{1 \over 
\kappa}} \over (\log{1 \over \kappa})^{\lambda}} \Bigr)$$
because $\int_{0}^{\infty} ({\sin{v} \over v})^2 \log{v} dv = -{\pi \over 2}
(C_0 + \log{2} - 1)$ (see [\ref{GR}] page $590$). The lemma follows from the
relationship between $C'$ and $D$.
\begin{lemma}
\label{lemma2.3}
If $K$ is even, $K''$ continuous, $\int_{-\infty}^{\infty} |K| < \infty$, $K(x)
\rightarrow 0$ as $x \rightarrow +\infty$, $K' \rightarrow 0$ as $x \rightarrow
+\infty$, and if $K''(x) \ll x^{-3}$ as $x \rightarrow +\infty$, then
\begin{equation}
\label{l3}
{\hat K}(t) = \int_{0}^{\infty} K''(x)\Bigl({\sin{\pi t x}\over \pi t}\Bigr)^2
dx.
\end{equation}
\end{lemma}

Proof: This is Lemma 3 of Goldston and Montgomery [\ref{GM}]. One can justify
it by integration by parts twice.

In Goldston and Montgomery [\ref{GM}], they used the function
\begin{equation}
\label{kernel1}
K_{\eta}(x) = {\sin{2\pi x} + \sin{2\pi (1+\eta)x} \over 2\pi x (1-4{\eta}^2
x^2)}
\end{equation}
for $\eta > 0$. Then
\begin{equation}
\label{kernel2}
{\hat K}_{\eta}(t)= \left\{ \begin{array}{ll}
 1, & \mbox{if $|t| \leq 1$},\\
 \cos^2{({\pi(|t|-1) \over 2\eta})}, & \mbox{if $1 \leq |t| \leq 1+\eta$},\\
 0, & \mbox{if $|t| \geq 1+\eta$}.
\end{array} \right.
\end{equation}

\begin{lemma}
\label{lemma2.35}
For $0 < \eta < {1 \over 4}$ and $K_{\eta}(x)$ defined as in (\ref{kernel1}),
then $K_{\eta}(x)$ satisfies the conditions in Lemma \ref{lemma2.3}. In
particular,
$${K_\eta}'' (x) \ll min\Bigl(1, {1 \over \eta^3 x^3}\Bigr).$$
Here the implicit constant is absolute.
\end{lemma}

Proof: One can easily see that $K_\eta$ is even, $\int_{-\infty}^{\infty}
|K_\eta| < \infty$, and $K_\eta \rightarrow 0$ as $x \rightarrow +\infty$.
Breaking it down into partial fractions, we have
\begin{equation*}
\begin{split}
K_\eta =& [\sin{2\pi x} + \sin{2\pi(1+\eta)x}] \Bigl[{1 \over x} + {\eta \over
1 - 2\eta x} + {-\eta \over 1 + 2\eta x}\Bigr], \\
{K_\eta}' =& 2\pi [\cos{2\pi x} + (1+\eta) \cos{2\pi(1+\eta)x}] \Bigl[{1 \over
x} + {\eta \over 1 - 2\eta x} + {-\eta \over 1 + 2\eta x}\Bigr] \\
&+ [\sin{2\pi x} + \sin{2\pi(1+\eta)x}] \Bigl[{-1 \over x^2} + {2\eta^2 \over
(1-2\eta x)^2} + {2\eta^2 \over (1+2\eta x)^2} \Bigr], \\
{K_\eta}'' =& 4\pi^2 [-\sin{2\pi x} - (1+\eta)^2 \sin{2\pi(1+\eta)x}] \Bigl[{1
\over x} + {\eta \over 1 - 2\eta x} + {-\eta \over 1 + 2\eta x}\Bigr] \\
&+ 4\pi [\cos{2\pi x} + (1+\eta) \cos{2\pi(1+\eta)x}] \Bigl[{-1 \over x^2} + {2
\eta^2 \over (1-2\eta x)^2} + {2\eta^2 \over (1+2\eta x)^2} \Bigr] \\
&+ [\sin{2\pi x} + \sin{2\pi(1+\eta)x}] \Bigl[{2 \over x^3} + {8\eta^3 \over
(1-2\eta x)^3} + {-8\eta^3 \over (1+2\eta x)^3} \Bigr].
\end{split}
\end{equation*}
Clearly, ${K_\eta}' \rightarrow 0$ as $x \rightarrow +\infty$. Moreover,
${K_\eta}''$ is continuous everywhere except possibly the points $0$, ${1 \over
2\eta}$ and $-{1\ \over 2\eta}$. Now, using Taylor expansions for $\sin$ and
$\cos$, we have, for $-1 \leq x \leq 1$,
\begin{eqnarray*}
{K_\eta}'' &=& 4\pi^2 [-2\pi x - 2\pi(1+\eta)^3 x + O(x^3)] \Bigl[{1 \over x} +
{\eta \over 1 - 2\eta x} + {-\eta \over 1 + 2\eta x} \Bigr] \\
& & + 4\pi [2 + \eta + O(x^2)] \Bigl[{-1 \over x^2} + {2 \eta^2 \over (1-2\eta x
)^2} + {2\eta^2 \over (1+2\eta x)^2} \Bigr] \\
& & + [2\pi(2+\eta)x + O(x^3)] \Bigl[{2 \over x^3} + {8\eta^3 \over
(1-2\eta x)^3} + {-8\eta^3 \over (1+2\eta x)^3} \Bigr]
\end{eqnarray*}
which shows that ${K_\eta}''$ is well-defined and continuous at $x=0$.
As for ${1 \over 2\eta}-1 \leq x \leq {1 \over 2\eta}+1$, we have
\begin{eqnarray*}
\sin{2\pi x} = \sin{(2\pi(x-{1 \over 2\eta}) + {\pi \over \eta})}, & &
\sin{2\pi(1+\eta)x} = -\sin{(2\pi(1+\eta)(x-{1 \over 2\eta}) + {\pi \over
\eta})}, \\
\cos{2\pi x} = \cos{(2\pi(x-{1 \over 2\eta}) + {\pi \over \eta})}, & &
\cos{2\pi(1+\eta)x} = -\cos{(2\pi(1+\eta)(x-{1 \over 2\eta}) + {\pi \over
\eta})}.
\end{eqnarray*}
So,
\begin{eqnarray*}
{K_\eta}'' &=& 4\pi^2 \Bigl[ -(\sin{\pi \over \eta})(1-(1+\eta)^2) -2\pi(
\cos{\pi \over \eta})(1 - (1+\eta)^3) (x-{1 \over 2\eta}) \\
& &+ O\Bigl((x-{1 \over 2\eta})^2\Bigr) \Bigr] \Bigl[{-{1 \over 2} \over x - {1
\over 2\eta}} + O(1) \Bigr] + 4\pi \Bigl[-\eta \cos{\pi \over \eta} -2\pi (\sin{
\pi \over \eta}) (1-(1+\eta)^2) \\
& &(x-{1 \over 2\eta}) - {(2\pi)^2 \over 2}(\cos{\pi \over \eta})
(1 - (1+\eta)^3) (x-{1 \over 2\eta})^2 + O\Bigl((x-{1 \over 2\eta})^3 \Bigr)
\Bigr] \Bigl[{{1 \over 2} \over (x-{1 \over 2\eta})^2} \\
& &+ O(1) \Bigr] + \Bigl[-2\pi \eta (\cos{\pi \over \eta})(x-{1 \over 2\eta}) -
{(2\pi)^2 \over 2} (\sin{\pi \over \eta}) (1-(1+\eta)^2) (x-{1 \over 2\eta})^2
\\
& &-{(2\pi)^3 \over 6} (\cos{\pi \over \eta}) (1-(1+\eta)^3) (x-{1 \over 2\eta})
^3 + O\Bigl((x-{1 \over 2\eta})^4\Bigr) \Bigr] \Bigl[{-1 \over (x-{1 \over 2
\eta})^3} + O(1) \Bigr].
\end{eqnarray*}
The ``pole terms'' cancel out exactly. So, we see that ${K_\eta}''$ is
well-defined and continuous at $x = {1 \over 2\eta}$. Replacing $\eta$ by
$-\eta$, one gets the well-definedness and continuity at $x = -{1 \over
2\eta}$.

The above calculations show that ${K_\eta}'' \ll 1$ when $0 \leq x \leq 1$ and
${1 \over 2\eta}-1 \leq x \leq {1 \over 2\eta}+1$. One can easily see that 
${K_\eta}'' \ll 1$ when $1 \leq x \leq {1 \over 2\eta}-1$. Hence,
\begin{equation}
\label{2.35.1}
{K_\eta}''(x) \ll 1 \ll {1 \over \eta^3 x^3} \mbox{ when } 0 \leq
x \leq {1 \over 2\eta}+1.
\end{equation}
When ${1 \over 2\eta}+1 \leq x \leq {1 \over \eta}$, after recombining the
partial fractions, one has
\begin{equation}
\label{K1}
{K_\eta}''(x) \ll {1 \over x(4\eta^2 x^2 -1)} + {12\eta^2 x^2 -1 \over x^2
(4\eta^2 x^2 -1)^2} + {1 \over x^3} + {1 \over (x-{1 \over 2\eta})^3}.
\end{equation}
But $4\eta^2 x^2-1 \geq (2\eta(1+{1 \over 2\eta}))^2-1 = (2\eta+1)^2 -1 \geq
\eta$, $\eta x \leq 1$, and $12\eta^2 x^2-1 \geq 4\eta^2 x^2$.
Hence,
\begin{equation}
\label{K2}
{1 \over x(4\eta^2 x^2 -1)} \leq {1 \over \eta x} \leq {1 \over \eta^3 x^3}
\mbox{ and } {12\eta^2 x^2 -1 \over x^2 (4\eta^2 x^2 -1)^2} \ll {x^2 \eta^2
\over x^2 \eta^2} \leq {1 \over \eta^3 x^3}.
\end{equation}
Meanwhile, $(1+2\eta)(1-\eta) = 1+\eta-2\eta^2 \geq 1$, so $x \geq {1 + 2\eta
\over 2\eta} \geq {1 \over 2\eta (1-\eta)}$. This implies that $x - {1 \over
2\eta} \geq \eta x$. Thus,
\begin{equation}
\label{K3}
{1 \over (x-{1 \over 2\eta})^3} \ll {1 \over \eta^3 x^3}.
\end{equation}
Combining (\ref{K1}), (\ref{K2}) and (\ref{K3}), we have
\begin{equation}
\label{2.35.2}
{K_\eta}''(x) \ll {1 \over \eta^3 x^3} \mbox{ when } {1 \over 2\eta}+1 \leq x
\leq {1 \over \eta}.
\end{equation}
Finally, when ${1 \over \eta} \leq x$, $4\eta^2 x^2 -1 \geq 2\eta^2 x^2$ and $x
- {1 \over 2\eta} \geq {x \over 2}$. Putting these into (\ref{K1}) (which is
still valid), we get
\begin{equation}
\label{2.35.3}
{K_\eta}''(x) \ll {1 \over \eta^2 x^3} + {\eta^2 x^2 \over x^2 \eta^4 x^4} + {1
\over x^3} \ll {1 \over \eta^3 x^3}.
\end{equation}
Consequently, the lemma follows from (\ref{2.35.1}), (\ref{2.35.2}) and
(\ref{2.35.3}).
\begin{lemma}
\label{lemma2.4}
If $f$ is a non-negative function defined on $[0,+\infty)$, $f(t) \ll 
{\log{(t+2)}}^2$, and if 
$$I(\kappa) = \int_{0}^{\infty} \Bigl({\sin{\kappa t} \over t} \Bigr)^2
f(t) dt = {\pi \over 2}\kappa \log{1 \over \kappa} + \pi C' \kappa + 
\epsilon(\kappa) \kappa$$ where $\epsilon(\kappa) = O \Bigl({1 \over (\log{1 /
\kappa})^5} \Bigr)$ for $T^{-1}(\log{T})^{-2} \leq \kappa \leq T^{-1}
(\log{T})^9$, then
$$J(T) = \int_{0}^{T} f(t) dt = T\log{T} + DT + O \Bigl({T \over \log{T}}
\Bigr).$$
Again, $C' = {1 \over 2}D - {1 \over 2}(C_0+\log{2})+1$.
\end{lemma}

Proof: Essentially, we follow Lemma 4 of Goldston and Montgomery [\ref{GM}].
Let $\eta = {1 \over (\log{T})^2}$. Set $X_1 = (\log{T})^{-2}$ and $X_2 = {1
\over 4}(\log{T})^9$. Let $K_{\eta}$ be the kernel defined as in
(\ref{kernel1}). Replace $t$ by $t/T$ in (\ref{l3}), multiply by $f(t) - \log{t}
-(D+1)$, and integrate over $0 \leq t < \infty$, one finds that
$$\int_{0}^{\infty} \bigl(f(t)-\log{t}-(D+1)\bigr){\hat K}_{\eta} \Bigl({t \over
T} \Bigr) dt = {T^2 \over {\pi}^2} \int_{0}^{\infty} {K_\eta}''(x) R\Bigl(
{\pi x \over T}\Bigr) dx$$
where
\begin{eqnarray*}
R(\kappa) &=& I(\kappa) - \int_{0}^{\infty} \Bigl({\sin{\kappa t}\over t}\Bigr)
^2 \log{t} dt
- (D+1) \int_{0}^{\infty} \Bigl({\sin{\kappa t} \over t}\Bigr)^2 dt \\
&=& I(\kappa) - \kappa \int_{0}^{\infty} \Bigl({\sin{\kappa t} \over \kappa
t}\Bigr)^2 
\log{\kappa t} d\kappa t - \kappa \log{1 \over \kappa} \int_{0}^{\infty} 
\Bigl({\sin{\kappa t} \over \kappa t}\Bigr)^2 d\kappa t \\
& &-(D+1) \int_{0}^{\infty} \Bigl({\sin{\kappa t} \over t}\Bigr)^2 dt\\
&=& I(\kappa) - {\pi \over 2} \kappa \log{1 \over \kappa}- {\pi \over 2} \kappa
\Bigl((D+1) - (C_0 + \log{2} -1) \Bigr),
\end{eqnarray*}
because $\int_{0}^{\infty} ({\sin{v} \over v})^2 dv = {\pi \over 2}$ and 
$\int_{0}^{\infty} ({\sin{v} \over v})^2 \log{v} dv= -{\pi \over 2}(C_0 +
\log{2}-1)$. Since,
$$I(\kappa) \ll \int_{0}^{\infty} min({\kappa}^2,t^{-2})(\log{(t+2)})^2 dt \ll
\kappa \Bigl(\log{(2+{1 \over \kappa})} \Bigr)^2$$
for all $\kappa > 0$, we see that, by Lemma \ref{lemma2.35},
${K_\eta}''(x)$ is bounded when $0 \leq x \leq X_1$. Thus,
$$\int_{0}^{X_1} {K_\eta}''(x)R \Bigl({\pi x \over T}\Bigr) dx \ll \int_{0}^
{X_1} {x \over T} \Bigl(\log{T \over x}\Bigr)^2 dx \ll {(\log{T})^2 {X_1}^2
\over T} \ll {1 \over T(\log{T})^2}.$$
On the other hand, by Lemma \ref{lemma2.35},
$$\int_{X_2}^{\infty} {K_\eta}''(x)R \Bigl({\pi x \over T}\Bigr) dx \ll
\int_{X_2}^{\infty} {1 \over \eta^3 x^3} {x \over T} (\log{T})^2 dx \ll
{(\log{T})^2 \over \eta^3 T X_2} \ll {1 \over T \log{T}}.$$
Finally,
\begin{eqnarray*}
\int_{X_1}^{X_2} {K_\eta}''(x)R \Bigl({\pi x \over T} \Bigr) dx &\ll&
\int_{X_1}^{X_2} |{K_\eta}''(x)|
\Big|\epsilon({\pi x \over T}) \Big| {\pi x \over T} dx \\ &\ll& 
\int_{X_1}^{1/\eta} {1 \over (\log{T})^5}{x \over T} dx + \int_{1/\eta}^{X_2}
{1 \over (\log{T})^5}{1 \over \eta^3 T x^2} dx \\
&\ll& {1 \over T\log{T}}
\end{eqnarray*}
by Lemma \ref{lemma2.35}. Therefore,
\begin{eqnarray*}
\int_{0}^{\infty} f(t) {\hat K}_{\eta} \Bigl({t \over T}\Bigr) dt &=&
\int_{0}^{\infty} \bigl(\log{t} +(D+1)\bigr){\hat K}_{\eta}\Bigl({t \over T}
\Bigr) dt +  O\Bigl({T \over \log{T}}\Bigr) \\
&=&\int_{0}^{T} \log{t}+(D+1) dt + O \Bigl(\int_{T}^{(1+\eta)T} \log{t}+(D+1) 
dt \Bigr)\\
& & + O\Bigl({T \over \log{T}}\Bigr) \\
&=& T\log{T} + DT + O(\eta T\log{T}) + O\Bigl({T \over \log{T}}\Bigr).
\end{eqnarray*}
Since $f$ is non-negative, we see that
$$\int_{0}^{\infty} f(t) {\hat K}_{\eta} \Bigl({(1+\eta)t \over T}\Bigr) dt \leq
J(T) \leq
\int_{0}^{\infty} f(t) {\hat K}_{\eta} \Bigl({t \over T} \Bigr) dt,$$
which gives the desired result.

\begin{lemma}
\label{lemma2.5}
For $0< \delta <1$, let
\begin{equation}
\label{a_s}
a(s) = {(1+\delta)^s -1 \over s}.
\end{equation}
If $|c(\gamma)| \leq 1$ for all $\gamma$ then
\begin{equation}
\label{l5}
\begin{split}
\int_{-\infty}^{+\infty} |a(it)|^2 \Big| \sum_{\gamma} {c(\gamma) \over
1+(t-\gamma)^2} \Big|^2 dt &= \int_{-\infty}^{+\infty} \Big|\sum_{|\gamma| \leq 
Z} {a({1 \over 2} + i\gamma) c(\gamma) \over 1+(t-\gamma)^2} \Big|^2 dt \\
&+ O\Bigl({\delta}^2 (\log{2 \over \delta})^3\Bigr) + O\Bigl({(\log{Z})^3 \over
Z} \Bigr)
\end{split}
\end{equation}
provided that $Z \geq {1 \over \delta}$. Here, the summations are over the
imaginary parts of the non-trivial zeros of the Riemann zeta function.
\end{lemma}

Proof: This is just Lemma 10 of Goldston and Montgomery [\ref{GM}].
\bigskip

We also need the explicit formula for $\psi(x)$ (see Davenport [\ref{D}] $\S
17$, as well as Goldston and Montgomery [\ref{GM}]) to get:
\begin{equation}
\label{explicit}
\begin{split}
\psi(x +\delta x) - \psi(x) - \delta x
&= - \sum_{|\gamma| \leq Z} a(\rho) x^{\rho} + O\Bigl((\log{x})min(1,{x \over Z 
\Vert x \Vert}) \Bigr) \\
&+ O \Bigl((\log{x})min(1,{x \over Z \Vert (1+\delta)x \Vert}) \Bigr)+
O \Bigl({x (\log{xZ})^2 \over Z} \Bigr)
\end{split}
\end{equation}
where $a(s)$ is as in (\ref{a_s}), and $\Vert \theta \Vert =
\mathop{min}_{n} |\theta - n|$ is the distance from $\theta$ to the nearest
integer. The error terms contribute a negligible amount if we take
$Z=X(\log{X})^2$. $\rho = \sigma + i\gamma$ denotes zeros of Riemann zeta
function.

\begin{lemma}
\label{lemma2.6}
For any $t \geq 0$,
$$ \sum_{\gamma} {1 \over 1 + (t-\gamma)^2} \ll \log{(t+2)} $$
where the sum is over the imaginary parts of the zeros of the Riemann zeta
function.
\end{lemma}

Proof: This is a lemma in Davenport [\ref{D}] on page 98.
\begin{lemma}
\label{lemma2.7}
Assume Riemann Hypothesis.
\begin{equation*}
\begin{split}
\int_{X}^{2X} \Big|\sum_{|\gamma| \leq Z} a(\rho)x^{\rho}\Big|^2 dx &=
\int_{X}^{2X} \bigl(\psi((1+\delta)x) - \psi(x) - \delta x \bigr)^2 dx \\
&+ O\Bigl(X^{3/2}\delta^{1/2}\log{2 \over \delta}\log{X} + X(\log{X})^2
\Bigr)
\end{split}
\end{equation*}
where $a(s)$ is defined by (\ref{a_s}) and $Z \geq X(\log{X})^2$.
\end{lemma}

Proof: From (\ref{explicit}), we have for $X \leq x \leq 2X$,
$$-\sum_{|\gamma| \leq Z} a(\rho)x^{\rho} = \psi\bigl((1+\delta)x \bigr) - 
\psi(x) - \delta x + O(\log{x}).$$
So
$$\int_{X}^{2X} \Big|\sum_{|\gamma| \leq Z} a(\rho)x^{\rho}\Big|^2 dx =
\int_{X}^{2X} \bigl(\psi((1+\delta)x) - \psi(x) - \delta x \bigr)^2 dx + 
\mbox{error}$$
where
\begin{eqnarray*}
\mbox{error} &=& O\Bigl(\int_{X}^{2X} |\psi((1+\delta)x) - \psi(x) - \delta x|
\log{x} dx\Bigr) +
O\Bigl(\int_{X}^{2X} (\log{x})^2 dx\Bigr) \\
&\ll& \Bigl(\int_{X}^{2X} (\psi((1+\delta)x) - \psi(x) - \delta x)^2 dx\Bigr)^
{1 \over 2} 
\Bigl(\int_{X}^{2X} (\log{x})^2 dx \Bigr)^{1 \over 2} + X(\log{X})^2 \\
&\ll& \Bigl(\delta X^2 (\log{2 \over \delta})^2 \Bigr)^{1 \over 2} X^{1 \over 2
} \log{X} + X(\log{X})^2 \\
&\ll& X^{3 \over 2} \delta^{1 \over 2} \log{2 \over \delta} \log{X} + X
(\log{X})^2
\end{eqnarray*}
by Cauchy's inequality and (\ref{ineq1}).
\section{Proof of (b) ``$\Leftrightarrow$'' (c)}
\begin{thm}
\label{thm2.1}
Assume Riemann Hypothesis. For $0< B_1 \leq B_2 < 1$, if
\begin{equation}
\label{2.1}
F(X,T) = {T \over 2 \pi}\log{T} + {D \over 2\pi}T + O\Bigl({T \over
(\log{T})^8}\Bigr)
\end{equation}
holds uniformly for $X^{B_1}(\log{X})^{-11} \leq T \leq X^{B_2}(\log{X})^{11}$,
then
\begin{equation}
\label{2.2}
\int_{1}^{X} \bigl(\psi(x+\delta x)-\psi(x)-\delta x \bigr)^2 dx =
   {1 \over 2}\delta X^2 \log{1 \over \delta} + C\delta X^2 + 
   O \Bigl({\delta X^2 \over \log{1/ \delta}} \Bigr)
\end{equation}
holds uniformly for $X^{-B_2} \leq \delta \leq X^{-B_1}$.

Conversely, for $1 < A_1 \leq A_2 < \infty$, if 
\begin{equation}
\label{2.3}
\int_{1}^{X} \bigl(\psi(x+\delta x)-\psi(x)-\delta x \bigr)^2 dx =
   {1 \over 2}\delta X^2 \log{1 \over \delta} + C\delta X^2 + 
   O \Bigl({\delta X^2 \over (\log{1/ \delta})^5} \Bigr)
\end{equation}
holds uniformly for $X^{-1/A_1}(\log{X})^{-10} \leq \delta \leq
X^{-1/A_2}(\log{X})^{10}$, then
\begin{equation}
\label{2.4}
F(X,T) = {T \over 2 \pi}\log{T}+ {D \over 2\pi}T + O \Bigl({T \over \log{T}}
\Bigr)
\end{equation}
holds uniformly for $T^{A_1} \leq X \leq T^{A_2}$.
Here $C$ and $D$ are related by $C = {D \over 2} - {C_0 \over 2} + 1$.
\end{thm}

Proof: We follow Goldston and Montgomery [\ref{GM}] closely. First we assume
(\ref{2.1}) and derive (\ref{2.2}). Let
$$J(T) = J(X,T) = 4 \int_{0}^{T} \Big|\sum_{\gamma} {X^{i\gamma} \over
1+(t-\gamma)^2} \Big|^2 dt.$$
Montgomery [\ref{M}] (see his (26), but be aware of the changes in notation)
showed that
$$J(X,T) = 2 \pi F(X,T) + O\bigl((\log{T})^3\bigr).$$
Thus, (\ref{2.1}) is equivalent to
\begin{equation}
\label{2.0.1}
J(X,T) = T\log{T} + DT + O\Bigl({T \over (\log{T})^8}\Bigr).
\end{equation}
With $a(s)$ defined by (\ref{a_s}), we note that
$$|a(it)|^2 = 4 \Bigl({\sin{\kappa t} \over t}\Bigr)^2$$
where $\kappa = {1 \over 2}\log{(1+\delta)}$. We have $\kappa = {\delta \over
 2} + O(\delta^2)$ and $\log{1 \over \kappa} = \log{1 \over \delta} + \log{2} +
O(\delta)$. Then, by Lemma \ref{lemma2.2} with $\lambda = 8$ and Lemma
\ref{lemma2.6}, we deduce that
\begin{equation}
\label{2.0.2}
\begin{split}
& \int_{0}^{\infty} |a(it)|^2 \Big|\sum_{\gamma} {X^{i\gamma} \over 1+(t-\gamma)
^2} \Big|^2 dt = {\pi \over 2} \kappa \log{1 \over \kappa} + \pi C' \kappa +
O\Bigl({\kappa \over (\log{1 / \kappa})^7} \Bigr) \\
&= {\pi \over 2}{\delta \over 2}\Bigl(\log{1 \over \delta} +\log{2} +O(\delta)
\Bigr) + \pi C' {\delta \over 2} + O(\delta^2) + O\Bigl({\delta \over
(\log{1/\delta})^7}\Bigr) \\
&= {\pi \over 4}\delta \log{1 \over \delta} +\Bigl({\pi \log{2} \over 4} + 
{\pi C' \over 2}\Bigr)\delta + O\Bigl({\delta \over (\log{1/\delta})^7}\Bigr).
\end{split}
\end{equation}
The values of $T$ for which we have used (\ref{2.1}) lie in the range
\begin{equation}
\label{cond1}
\delta^{-1}\Bigl(\log{1 \over \delta}\Bigr)^{-10} \leq T \leq 3\delta^{-1} 
\Bigl(\log{1 \over \delta} \Bigr)^{10}.
\end{equation}
The integrand is even, so the value is doubled if we integrate over
negative values of $t$ as well. Then, by Lemma \ref{lemma2.5},
$$\int_{-\infty}^{+\infty} \Big|\sum_{|\gamma| \leq Z} {a(\rho) X^{i \gamma} 
\over
1+(t-\gamma)^2} \Big|^2 dt = {\pi \over 2}\delta \log{1 \over \delta} + 
\Bigl({\pi \log{2} \over 2}+ \pi C'\Bigr) \delta + O\Bigl({\delta \over 
(\log{1/\delta})^7}\Bigr)$$
provided that $Z \geq \delta^{-1}(\log{1 \over \delta})^{10}$. Let $S(t)$
denote the above sum over $\gamma$. Its Fourier transform is
$${\hat S}(u) = \int_{-\infty}^{+\infty} S(t)e(-tu)dt = \pi \sum_{|\gamma| \leq
Z} a(\rho) X^{i\gamma} e(-\gamma u) e^{-2 \pi |u|}.$$
Hence, by Plancherel's identity, the above integral is
$$= \pi^2 \int_{-\infty}^{+\infty} \Big|\sum_{|\gamma| \leq Z} a(\rho)
X^{i\gamma}\Big|^2 e^{-4\pi |u|} du.$$
Let $Y=\log{X}$, $-2\pi u = y$, we have
\begin{equation}
\label{2.0.3}
\int_{-\infty}^{+\infty} \Big|\sum_{|\gamma| \leq Z} a(\rho) e^{i\gamma (Y+y)}
\Big|^2
e^{-2|y|} dy = \delta \log{1 \over \delta}+(\log{2} + 2C') \delta + 
O\Bigl({\delta \over (\log{1/\delta})^7}\Bigr).
\end{equation}
Now, set
$$f(y) = {|\sum_{|\gamma| \leq Z} a(\rho) e^{i\gamma y}|^2 \over
\delta \log{1 \over \delta} + (\log{2} + 2C') \delta},$$
then $f(y) \geq 0$ if $\delta$ is small enough. From (\ref{2.0.3}), we know
that $\int_{0}^{1} f(Y+y) dy \ll 1$ as well as $\int_{-\infty}^{+\infty} 
e^{-2|y|} f(Y+y) dy = 1 + O({1 \over Y^8})$ for $\log{1 \over \delta} \asymp 
\log{X} = Y$.
Therefore we can use Lemma \ref{lemma2.1}. With the change of variable $x =
e^{Y+y}$, we get
\begin{eqnarray*}
\int_{X}^{2X} \Big|\sum_{|\gamma| \leq Z}a(\rho) x^{\rho}\Big|^2 dx &=& 
\Bigl({3 \over 2} + O({1 \over (\log{1/\delta})^2}) \Bigr) \Bigl(\delta 
\log{1 \over \delta} + (\log{2}+2C')\delta \\
& &+ O\bigl({\delta \over (\log{1/\delta})^7}\bigr) \Bigr) X^2 \\
&=& {3 \over 2}\delta X^2 \log{1 \over \delta} + {3 \over 2}(\log{2}+2C')\delta
X^2 + O \Bigl({\delta X^2 \over \log{1/\delta}}\Bigr).
\end{eqnarray*}
By Lemma \ref{lemma2.7} with $Z=X(\log{X})^2$,
\begin{eqnarray*}
& &\int_{X}^{2X} \bigl(\psi(x+\delta x)-\psi(x)-\delta x \bigr)^2 dx \\
&=&{3 \over 2}\delta X^2 \log{1 \over \delta} + {3 \over 2}(\log{2}+2C')
\delta X^2 \\
& &+ O\Bigl({\delta X^2 \over \log{1/\delta}}\Bigr) + 
O\Bigl(X^{3/2}\delta^{1/2}\log{2 \over \delta}\log{X} + X(\log{X})^2 \Bigr).
\end{eqnarray*}
But as $X^{-B_2} \leq \delta \leq X^{-B_1}$, $0<B_1 \leq B_2 <1$, the first
error term is bigger than the other one.

Now, we replace $X$ by $X2^{-k}$ and sum over $1 \leq k \leq
K=[{2\log{\log{X}} \over \log{2}}]$. One obtains
\begin{eqnarray*}
\int_{X2^{-K}}^{X} \bigl(\psi(x+\delta x) - \psi(x) - \delta x\bigr)^2 dx &=&
{1 - 2^{-2K} \over 2} \Bigl(\delta X^2 \log{1 \over \delta} + (\log{2}+2C')
\delta X^2 \Bigr) \\
& &+ O\Bigl({\delta X^2 \over \log{1/\delta}}\Bigr).
\end{eqnarray*}
To bound the contribution from the range $1 \leq x \leq X/2^{K}$, we use
inequality (\ref{ineq1}):
$$\int_{1}^{X/2^{K}} \bigl(\psi(x+\delta x) - \psi(x) - \delta x \bigr)^2 dx 
\ll \delta {X^2 \over (\log{X})^4} \Bigl(\log{2 \over \delta}\Bigr)^2 \ll 
{\delta X^2 \over \log{(1/\delta)}}.$$
We finally get (\ref{2.2}) as $C=C' + {\log{2} \over 2}$. The last thing to
take care of is the range for $T$ and $\delta$. Now for $X^{-B_2} \leq \delta
\leq X^{-B_1}$, putting this into (\ref{cond1}), the whole argument goes
through as long as (\ref{2.1}) is true for $X^{B_1}(\log{X})^{-10} \ll T \ll 
X^{B_2} (\log{X})^{10}$. But we assume that this is true for $X^{B_1}(\log{X})
^{-11} \leq T \leq X^{B_2}(\log{X})^{11}$ to start with! Therefore (\ref{2.1})
implies 
(\ref{2.2}) when $X$ is sufficiently large.

We now deduce (\ref{2.4}) from (\ref{2.3}). Let $X_1 = X$ and $X_2 = X(\log{X})
^4$. Then
\begin{equation}
\label{more}
\int_{1}^{\tilde{X}} \bigl(\psi(x+\delta x) - \psi(x) -\delta x\bigr)^2 dx =
{1 \over 2}\delta {\tilde{X}}^2 \log{1 \over \delta} + C \delta {\tilde{X}}^2 +
O\Bigl({\delta {\tilde{X}}^2 \over (\log{1/\delta})^5}\Bigr)
\end{equation}
for $X_1 \leq \tilde{X} \leq X_2$ whenever
\begin{equation}
\label{cond2}
X^{-1/A_1}(\log{X})^{-6} \leq \delta \leq X^{-1/A_2}(\log{X})^{10}.
\end{equation}
By integration by parts and (\ref{more}), we find that
\begin{equation*}
\begin{split}
\int_{X_1}^{X_2} \bigl(\psi(x+\delta x)-\psi(x)-\delta x\bigr)^2 x^{-4} dx =&
{1 \over 2}{X_1}^{-2} \delta \log{1 \over \delta} + C {X_1}^{-2}\delta \\
+& O\Bigl({\delta {X_1}^{-2} \over (\log{1/\delta})^5} + {\delta {X_1}^{-2}
\log{1/\delta} \over (\log{X_1})^8}\Bigr).
\end{split}
\end{equation*}
From (\ref{ineq1}), we deduce that
$$\int_{X_2}^{\infty} \bigl(\psi(x+\delta x)-\psi(x)-\delta x\bigr)^2 x^{-4} dx
\ll {\delta {X_1}^{-2} (\log{1/\delta})^2 \over (\log{X_1})^8}.$$
We add these relations, and multiply through by ${X_1}^2$. By making a further
appeal to (\ref{more}) with $\tilde{X}=X_1$, we deduce that
\begin{eqnarray*}
& &\int_{0}^{\infty} min\Bigl({x^2 \over {X_1}^2},{{X_1}^2 \over x^2}\Bigr)
\bigl(\psi(x+\delta x)-\psi(x)-\delta x\bigr)^2 x^{-2} dx \\
&=& \delta \log{1 \over \delta} + 2C \delta + O\Bigl({\delta \over
(\log{1/\delta})^5}\Bigr).
\end{eqnarray*}
Write $X$ for $X_1$, put $Y=\log{X}$, $x=e^{Y+y}$, and appeal to the
explicit formula (\ref{explicit}) with $Z=X(\log{X})^2$ (the error term arises
from the explicit formula can be treated similarly as Lemma \ref{lemma2.7} and 
it is small because of the range of $\delta$). We have
$$\int_{-\infty}^{+\infty} \Big|\sum_{|\gamma| \leq Z} a(\rho) e^{i\gamma (Y+y)
}\Big|^2
e^{-2|y|} dy = \delta \log{1 \over \delta} + 2C \delta + O\Bigl({\delta
\over (\log{1/\delta})^5}\Bigr)$$
which is almost (\ref{2.0.3}) except for the error term. Now retracing our
steps, we have (\ref{2.0.2}) except that the error term is 
$O({\kappa \over (\log{1/\kappa})^5})$. Because of Lemma \ref{lemma2.6}, we can
use Lemma \ref{lemma2.4} and obtain (\ref{2.0.1}) except that the error term is
$O({T \over \log{T}})$. This works provided that $T^{-1}(\log{T})^{-2} \leq 
\kappa \leq T^{-1}(\log{T})^9$. However, as $T^{A_1} \leq X \leq T^{A_2}$, from
(\ref{cond2}), we have $X^{-1/A_1}(\log{X})^{-6} \leq \delta \leq
X^{-1/A_2}(\log{X})^{10}$ which allows $\kappa$ to be in the range we want as
$\delta \asymp \kappa$! So, the second half of the theorem is true.
\bigskip

One can modify the lemmas to get analogous results with smaller error terms.
\begin{thm}
Assume Riemann Hypothesis. For $0< B_1 \leq B_2 < 1$, if, for some
$\epsilon>0$,
\begin{equation*}
F(X,T) = {T \over 2 \pi}\log{T} + {D \over 2\pi}T + O(T^{1-\epsilon})
\end{equation*}
holds uniformly for $X^{B_1-{\epsilon}_1} \leq T \leq X^{B_2+{\epsilon}_1}$, 
then
\begin{equation*}
\int_{1}^{X} \bigl(\psi(x+\delta x)-\psi(x)-\delta x\bigr)^2 dx =
   {1 \over 2}\delta X^2 \log{1 \over \delta} + C\delta X^2 + 
   O(\delta^{1+{\epsilon}_2} X^2)
\end{equation*}
holds uniformly for $X^{-B_2} \leq \delta \leq X^{-B_1}$ where ${\epsilon}_1,
{\epsilon}_2 >0$ may depend on $\epsilon$.

Conversely, for $1 < A_1 \leq A_2 < \infty$, if, for some $\epsilon >0$, 
\begin{equation*}
\int_{1}^{X} \bigl(\psi(x+\delta x)-\psi(x)-\delta x\bigr)^2 dx =
   {1 \over 2}\delta X^2 \log{1 \over \delta} + C\delta X^2 + 
   O(\delta^{1+\epsilon} X^2)
\end{equation*}
holds uniformly for $X^{-1/A_1-{\epsilon}_1} \leq \delta \leq
X^{-1/A_2+{\epsilon}_1}$, then
\begin{equation*}
F(X,T) = {T \over 2 \pi}\log{T} + {D \over 2\pi}T + O(T^{1-{\epsilon}_2})
\end{equation*}
holds uniformly for $T^{A_1} \leq X \leq T^{A_2}$ where ${\epsilon}_1,
{\epsilon}_2 >0$ may depend on $\epsilon$.
Here $C$ and $D$ are related by $C = {D \over 2} - {C_0 \over 2} +1$.
\end{thm}


\end{document}